\documentclass[12pt]{article}
\usepackage{graphicx}
\usepackage{amsmath}
\usepackage{amsfonts}
\usepackage{amssymb}
\newtheorem{theorem}{Theorem}

\newtheorem{corollary}[theorem]{Corollary}

\newtheorem{proposition}[theorem]{Proposition}

\newenvironment{proof}[1][Proof]{\textit{#1.} }

\def\text{\hbox} 

\def\a{\alpha}
\def\b{\beta}
\def\bb{\gamma}
\def\ba{\tilde\alpha}
\def\wa{\widehat\alpha}
\def\la{\bar\alpha}
\def\lb{\bar\beta}
\def\lbb{\bar\gamma}
\def\g{\xi}
\def\d{d}
\def\dd{\delta}

\def\ee{\varepsilon}

\def\p{\pi}

\def\r{\rho}

\def\ps{\varphi}

\def\f{\psi}

\def\o{\omega}

\def\F{\Phi}

\def\O{\Omega}
\def\S{{\bf S}}

\def\A{{\cal A}}
\def\B{{\cal B}}
\def\K{{\cal K}}
\def\A'{{\cal A}'}

\def\Z{{\mathbb Z}}

\begin{document}

\title{\textbf{Strongly-cyclic branched coverings of
$\mathbf{(1,1)}$-knots and cyclic presentations of groups}}

\author{ALESSIA CATTABRIGA \textsc{and }MICHELE MULAZZANI\\\ \\\textit{Department of
Mathematics, University of Bologna,}\\ \textit{P.za di P.ta
S.Donato} 5, I-40126 \textit{Bologna, Italy.}\\ \textit{e-mail:}
\textsf{cattabri@dm.unibo.it  mulazza@dm.unibo.it}}

\maketitle

\begin{abstract}
{We study the connections among the mapping class group of the
twice punctured torus, the cyclic branched coverings of
$(1,1)$-knots and the cyclic presentations of groups. We give the
necessary and sufficient conditions for the existence and
uniqueness of the $n$-fold strongly-cyclic branched coverings of
$(1,1)$-knots, through the elements of the mapping class group. We
prove that every $n$-fold strongly-cyclic branched covering of a
$(1,1)$-knot admits a cyclic presentation for the fundamental
group, arising from a Heegaard splitting of genus $n$. Moreover, we
give an algorithm to produce the cyclic presentation and illustrate
it in the case of cyclic branched coverings of
torus knots of type $(k,hk\pm 1)$.  
}

\end{abstract}

\section{Introduction and preliminaries}

The problem of determining whether a balanced presentation of a
group is geometric (i.e. induced by a Heegaard diagram of a closed
orientable 3-manifold) is of considerable interest in geometric
topology and has already been examined by many authors (see
\cite{Gr, Mo, Ne, OS1, OS2, OS3, St}). Moreover, the connections
between cyclic coverings of $\S^3$ branched over knots and cyclic
presentations of their fundamental groups, induced by suitable
Heegaard diagrams, have recently been discussed in several papers
(see \cite{BKM, CHK1, CHK2, CHR, Du, HKM1, HKM2, Ki, KKV1, KKV2,
MR, VK}).

We recall that a finite balanced presentation of a group
$<x_1,\ldots,x_n\mid r_1,\ldots,r_n>$ is said to be a {\it cyclic
presentation\/} if there exists a word $w$ in the free group $F_n$
generated by $x_1,\ldots,x_n$ such that the relators of the
presentation are $r_k=\theta_n^{k-1}(w)$, $k=1,\ldots,n$, where
$\theta_n :F_n\to F_n$ is the automorphism defined by $\theta_n
(x_i)=x_{i+1}$ (subscripts mod $n$), $i=1,\ldots,n$. This cyclic
presentation and the related group will be denoted by $G_n(w)$, so
$$G_n(w)=<x_1,\ldots,x_n\mid
w,\theta_n(w),\ldots,\theta_n^{n-1}(w)>.$$ Obviously $G_n(w)\cong
G_n(\theta_n^s(w))$ for every integer $s$. The {\it polynomial
associated to the cyclic presentation\/} $G_n (w)$ is defined by
$$ f_w (t)= \sum_{i=1}^{n} a_i t^{i-1}, $$ where $a_i$ is the
exponent sum of $x_i$ in $w$. For further details see \cite{Jo}.

\medskip

Some interesting examples of cyclically presented groups which are
the fundamental groups of cyclic branched coverings of $\S^3$ are the following:
\begin{itemize}
\item[(i)] the Fibonacci group
$F(2n)=G_{2n}(x_1x_2x_3^{-1})=G_n(x_1^{-1}x_2^2x_3^{-1}x_2)$ is
the fundamental group of the $n$-fold cyclic covering of $\S^3$
branched over the figure-eight knot, for all $n>1$ (see
\cite{HKM2});
\item[(ii)] the Sieradski group
$S(n)=G_n(x_1x_3x_2^{-1})$ is the fundamental group of the
$n$-fold cyclic covering of $\S^3$ branched over the trefoil knot,
for all $n>1$ (see \cite{CHK1});
\item[(iii)] the fractional Fibonacci group
$\widetilde F_{l,k}(n)=G_n((x_1^{-l}x_2^l)^kx_2(x_3^{-l}x_2^l)^k)$ is
the fundamental group of the $n$-fold cyclic covering of $\S^3$
branched over the genus one two-bridge knot with Conway
coefficients $[2l,-2k]$, for all $n>1$ and $l,k>0$ (see
\cite{VK}).
\end{itemize}
Notice that all the above cyclic presentations are geometric
(i.e., they arise from suitable Heegaard diagrams).

In order to investigate the relations between cyclic branched
coverings of knots in $\S^3$ and cyclic presentations for their
fundamental groups, Dunwoody introduced in \cite{Du} a class of
Heegaard diagrams depending on six integers, having cyclic
symmetry and defining cyclic presentations for the corresponding
fundamental groups. In \cite{GM} it has been shown that the
3-manifolds represented by these diagrams -- the so-called {\it
Dunwoody manifolds\/} -- are cyclic coverings of lens spaces,
branched over $(1,1)$-knots. As a corollary, it has been proved
that for some well-determined cases the manifolds turn out to be
cyclic coverings of $\S^3$, branched over knots. This gives a
positive answer to a conjecture formulated by Dunwoody in
\cite{Du}, which has also been independently proved in \cite{SK}.

In what follows, we shall deal with $(1,1)$-knots, also called
genus one 1-bridge knots. They are knots in lens spaces (possibly
in $\S^3$) admitting the following decomposition, called
\textit{(1,1)-decomposition}. A knot $K\subset L(p,q)$ is said to
be a \textit{(1,1)-knot} if there exists a Heegaard splitting of
genus one
$$(L(p,q),K)=(T,A)\cup_{\ps}(T',A'),$$ where $T$ and $T'$ are solid
tori, $A\subset T$ and $A'\subset T'$ are properly embedded
trivial arcs\footnote{This means that
there exists a disk $D\subset T$ (resp. $D'\subset T'$) with
$A\cap D=A\cap\partial D=A$ and $\partial D-A\subset\partial T$
(resp. $A'\cap D'=A'\cap\partial D'=A'$ and $\partial
D'-A'\subset\partial T'$).}, and $\ps:(\partial T',\partial A')\to(\partial
T,\partial A)$ is an orientation-reversing (attaching) homeomorphism.

\bigskip
\begin{figure}[ht]
\begin{center}
\includegraphics*[totalheight=3cm]{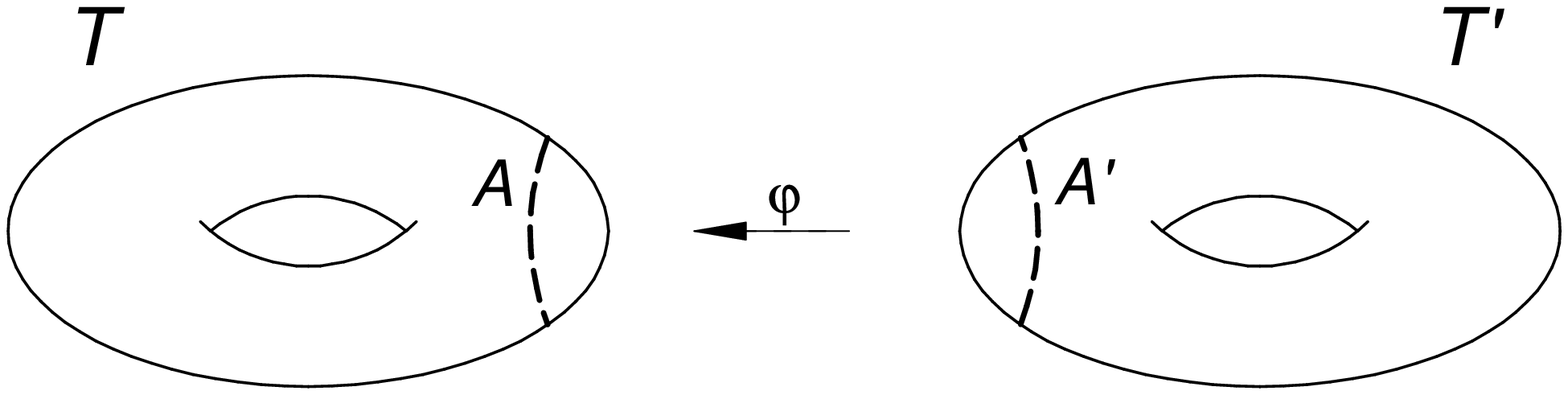}
\end{center}
\caption{A $(1,1)$-knot decomposition.}

\label{Fig. 1}

\end{figure}

Note that $(1,1)$-knots are just a particular case of the notion
of $(g,b)$-links in closed orientable 3-manifolds (see \cite{Do}
and \cite{Ha}), which generalize the classical concept of bridge
decomposition of links in $\S^3$. The class of $(1,1)$-knots is
very important in the light of some results and conjectures
involving Dehn surgery on knots (see references in \cite{Ha}). It is well known
that the subclass of $(1,1)$-knots in $\S^3$ contains all torus
knots (trivially) and all 2-bridge knots \cite{KS}.

In this paper we prove that every $n$-fold strongly-cyclic
branched covering of a $(1,1)$-knot admits a Heegaard splitting of
genus $n$ which encodes a cyclic presentation for the fundamental
group (see Section 4). The definition of strongly-cyclic branched
coverings of $(1,1)$-knots will be given in Section 3, where we
give necessary and sufficient conditions for their existence and
uniqueness, using a representation of $(1,1)$-knots by elements of
the mapping class group of the twice punctured torus (see Section
2). In Section 4 we give a constructive algorithm to obtain the
word defining the cyclic presentation. Moreover, we show that the
Alexander polynomial $\Delta_K(t)$ of a $(1,1)$-knot
$K\subset\S^3$ is equal to the polynomial associated to the cyclic
presentation of the fundamental group of the $n$-fold
strongly-cyclic branched covering of $K$, up to units of $\Z [t,
t^{-1}]$, when $n$ is greater than the degree of $\Delta_K(t)$.

\section{$\mathbf{(1,1)}$-knots and the mapping class group of the twice punctured torus}

\begin{figure}[ht]
\begin{center}
\includegraphics*[totalheight=4cm]{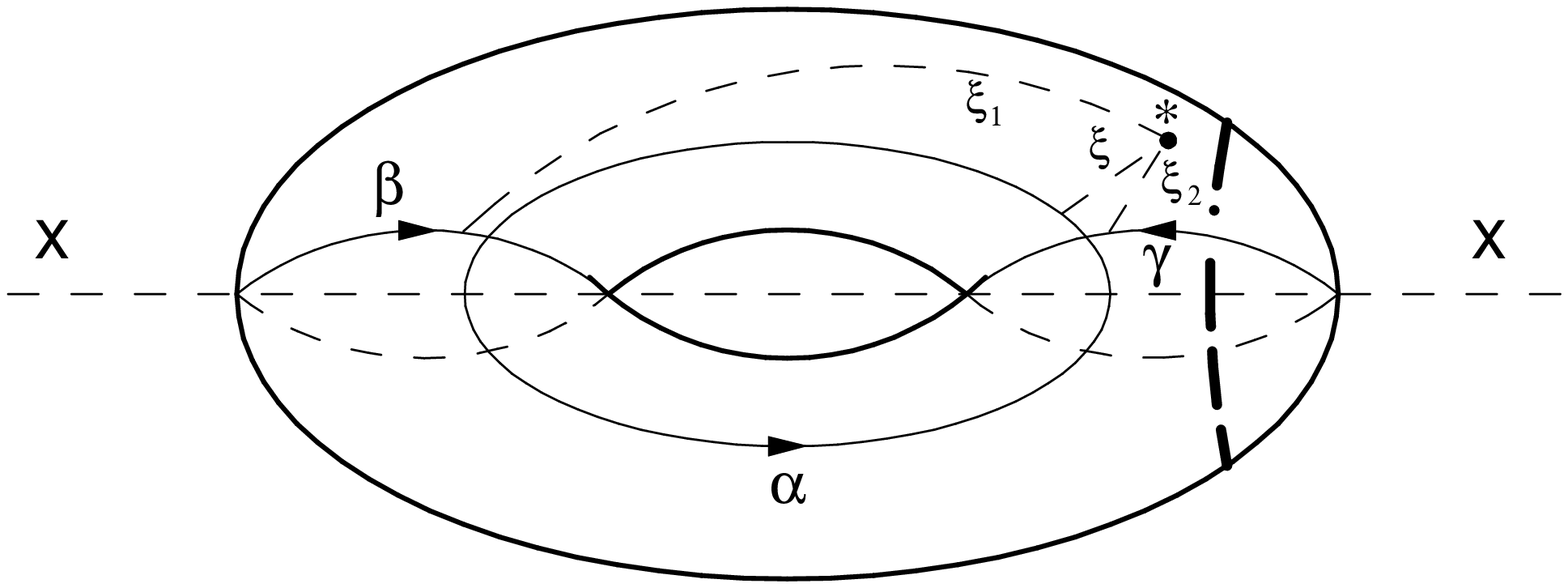}
\end{center}
\caption{Generators of $MCG(T_2)$.}

\label{Fig. 2}

\end{figure}

Let $K\subset L(p,q)$ be a $(1,1)$-knot with $(1,1)$-decomposition
$(L(p,q),K)=(T,A)\cup_{\ps}(T',A')$ and let $\tau :(T,A)\to
(T',A')$ be a fixed orientation-reversing homeomorphism, then
$\f=\ps\circ\tau_{|\partial T}$ is an orientation-preserving
homeomorphism of $(\partial T,\partial A)$. Moreover, since two
isotopic attaching homeomorphisms produce equivalent
$(1,1)$-knots, we have a natural surjective map from the
(orientation-preserving) mapping class group of the twice
punctured torus $MCG(T_2)$ to the class $\K_{1,1}$ of all
$(1,1)$-knots
$$\f\in MCG(T_2)\mapsto K_{\f}\in \K_{1,1}.$$ A standard set of
generators for $MCG(T_2)$ is given by the rotation $\r$ of $\p$
radians around the x--x axis and the three right-hand Dehn twists
$\d_{\a},\d_{\b},\d_{\bb}$, respectively around the curves
$\a,\b,\bb$, as depicted in Figure \ref{Fig. 2}. Observe that
$\d_{\a},\d_{\b},\d_{\bb}$ fix the punctures, while $\r$ exchanges
them.

The subgroup of $MCG(T_2)$ generated by
$\d_{\a},\d_{\b},\d_{\bb}$, i.e. the \textit{pure mapping class
group} of $T_2$, which contains only the elements fixing the
punctures is denoted by $PMCG(T_2)$ \cite{LP} (a very simple
presentation of $PMCG(T_2)$ can be found in \cite{PS}). Since $\r$
obviously commutes with the other generators, we have
$MCG(T_2)\cong PMCG(T_2)\oplus\Z_2$ and every element $\f$ of
$MCG(T_2)$ can be written as $\f=\f'\r^k$, where $\f'\in
PMCG(T_2)$. Since $\r$ can be extended to a homeomorphism of the
pair $(T,A)$, the $(1,1)$-knots $K_{\f}$ and $K_{\f'}$ are
equivalent. So, for our discussion it suffices to consider only
the elements of $PMCG(T_2)$.

The group $PMCG(T_2)$ naturally maps by an epimorphism to the
mapping class group of the torus (no punctures) $MCG(T_0)$, which
is generated by $\d_{\a}$ and $\d_{\b}=\d_{\bb}$. As is well
known, $MCG(T_0)$ is isomorphic to $SL(2,\Z)$ via the map which
associates to any element $\f\in MCG(T_0)$ the matrix representing
the isomorphism $\f_{\#}$ of $H_1(T_0)\cong\Z\oplus\Z$, with
respect to a fixed ordered basis $\B$. Assuming $\B=(\b,\a)$, we
have the epimorphism
\begin{equation}\label{isoMCG}
\O:PMCG(T_2)\to SL(2,\Z)
\end{equation}
given by $\d_{\b},\d_{\bb}\mapsto
\left(\begin{matrix}1&-1\\0&1\end{matrix}\right)$ and
$\d_{\a}\mapsto \left(\begin{matrix}1&0\\1&1\end{matrix}\right)$.
From this assumption, if
$\O(\f)=\left(\begin{matrix}q&s\\p&r\end{matrix}\right)$ then
$K_{\f}$ is a $(1,1)$-knot in the lens space $L(\vert p\vert,q)$ \cite{BZ}, and therefore
it is a knot in $\S^3$ if and only if $p=\pm 1$.

It is easy to see that if either $\f=1$ or $\f=\d_{\b}$ or $\f=\d_{\bb}$,
then $K_{\f}$ is the trivial knot in $\S^1\times\S^2$,
and if $\f=\d_{\a}$ then $K_{\f}$ is the trivial knot in $\S^3$.

Now we compute the fundamental group of the complement of a
$(1,1)$-knot, by applying the Seifert-Van Kampen theorem to its
$(1,1)$-decomposition. In order to do that, we fix the base point
$*$  and define the loops $\la=\g\cdot\a\cdot\g^{-1}$,
$\lb=\g_1\cdot\b\cdot\g_1^{-1}$ and
$\lbb=\g_2\cdot\bb\cdot\g_2^{-1}$, where $\g,\g_1,\g_2$ are the
paths connecting $*$ to $\a$, $\b$ and $\bb$ respectively, as
depicted in Figure \ref{Fig. 2}. It is easy to see that $\la,\lb$
and $\lbb$ are homologous to $\a,\b$ and $\bb$ respectively.
Moreover, we set $\la'=\tau(\la)$, $\lb'=\tau(\lb)$,
$\lbb'=\tau(\lbb)$ and $*'=\tau(*)$. Observe that
$*'=\ps^{-1}(*)$, up to isotopy. The homotopy classes of
$\la,\lb,\lbb$ generate $\p_1(\partial T-\partial A,*)$ and the
homotopy classes of $\la',\lb',\lbb'$ generate $\p_1(\partial
T'-\partial A',*')$. In order to simplify the notation, in the
following we use the same symbol for loops (resp. for cycles) and
for corresponding homotopy classes (resp. homology classes).

\begin{proposition} \label{fundamental} The fundamental group of
the complement of a $(1,1)$-knot
$K_{\f}\subset L(p,q)$ admits the presentation
$$\p_1(L(p,q)-K,*)=<\la,\lbb\,|\,r(\la,\lbb)>,$$ where
$r(\la,\lbb)$ is obtained by erasing all $\lb$'s terms from the
word representing the homotopy class of $\f(\lb)$.
\end{proposition}
\begin{proof}
 We have $\pi_1(T-A,*)=<\la,\lb,\lbb\,|\,\lb>$ and
$\pi_1(T'-A',*')=<\la',\lb',\lbb'\,|\,\lb'>$. Applying Seifert-Van
Kampen theorem we get
 \begin{eqnarray*}
 \p_1(L(p,q)-K,*)&=&<\la,\lb,\lbb,\la',\lb',\lbb'\,|\,
\lb,\lb',\ps(\la')\la^{\prime{-1}}, \ps(\lb')\lb^{\prime{-1}}, \ps(\lbb')\lbb^{\prime{-1}}>\\
 &=&<\la,\lb,\lbb\,|\,\lb,\f(\lb)>=<\la,\lbb\,|\,r(\la,\lbb)>,
 \end{eqnarray*}
and the proof is obtained.
\end{proof}

\medskip

From this result the computation of the first homology group is
straightforward, observing that $\{\a,\b,\bb\}$ is a set of free
generators for $H_1(\partial T-\partial A)$ and that, by
(\ref{isoMCG}), the homology relation $\f(\b)=p\a +q'\b +q''\bb$
(where $q'+q''=q$) holds for a $(1,1)$-knot $K_{\f}\subset
L(p,q)$.

\begin{corollary} \label{homology} Let $K_{\f}$ be a
$(1,1)$-knot in $L(p,q)$ then
$$H_1(L(p,q)-K_{\f})=<\a,\bb\,|\,p\a+q''\bb>\cong\Z\oplus\Z_{\gcd(p,q'')},$$
where $q''$ is uniquely determined by the homology relation:
$$\f(\b)=p\a+q'\b+q''\bb.$$
\end{corollary}

In order to compute the homotopy (and homology) classes of
$\f(\lb)$, we list the actions of the homomorphisms
$\d_{\a},\d_{\b},\d_{\bb}$ and their inverses on the homotopy
classes $\la,\lb,\lbb$.
\begin{equation}\label{formuleMCG}\begin{array}{lclcl}\d_{\a}(\la)=\la&\ &\d_{\a}(\lb)=\la\lb&\ &\d_{\a}(\lbb)=\la\lbb,\\
\d_{\a}^{-1}(\la)=\la&\ &\d_{\a}^{-1}(\lb)=\la^{-1}\lb&\ &\d_{\a}^{-1}(\lbb)=\la^{-1}\lbb,
\\\d_{\b}(\la)=\lb^{-1}\la&\ &\d_{\b}(\lb)=\lb&\ &\d_{\b}(\lbb)=\lbb,
\\\d_{\b}^{-1}(\la)=\lb\la&\ &\d_{\b}^{-1}(\lb)=\lb&\ &\d_{\b}^{-1}(\lbb)=\lbb,
\\\d_{\bb}(\la)=\la\lbb^{-1}&\ &\d_{\bb}(\lb)=\lb&\ &\d_{\bb}(\lbb)=\lbb,
\\\d_{\bb}^{-1}(\la)=\la\lbb&\ &\d_{\bb}^{-1}(\lb)=\lb&\ &\d_{\bb}^{-1}(\lbb)=\lbb.
\end{array}\end{equation}

\medskip

The next proposition describes a class of non-trivial
$(1,1)$-knots in $\S^3$.

\begin{proposition} \label{proposition 0}
If $\f=\d_{\a}^{\pm h}\d_{\bb}^{-k}\d_{\b}^{1+k}\d_{\a}$, with
$h,k\in\Z^+$, then $K_{\f}$ is the torus knot $\mathbf{t}(k,hk\mp 1)$
in $\S^3$. In particular, for $h=1$ and $k=2$, the knot $K_{\f}$ is the trefoil knot.
\end{proposition}
\begin{proof} $K_{\f}$ is always a knot in $\S^3$,
since $\O(\f)=\left(\begin{matrix}0&-1\\1&1\mp
h\end{matrix}\right)$. Then, using (\ref{formuleMCG}), we obtain
$\f(\lb)=(\lb^{-1}\la^{\mp h})^{1+k}\la(\la^{\pm h}\lbb)^k\la^{\pm
h}\lb$ and therefore $r(\la,\lbb)=\la^{1\mp hk}(\lbb\la^{\pm
h})^k$. We have
$$\pi_1(\S^3-K_{\f},*)=<\la,\lbb\,|\,\la^{1\mp hk}(\lbb\la^{\pm
h})^k>\cong <\la,u\,|\,\la^{1\mp hk}u^k>,$$
 which is isomorphic to
the group of the torus knot $\mathbf{t}(k,hk\mp 1)$. From a result
of \cite{BM}, we can conclude that in fact
$K_{\f}=\mathbf{t}(k,hk\mp 1)$.
\end{proof}

\section{Strongly-cyclic branched coverings of $\mathbf{(1,1)}$-knots}

An $n$-fold cyclic branched covering between two orientable closed
manifolds $f:M\to N$, with branching set $L$, is completely
defined by an epimorphism (called {\it monodromy\/})
$\o_{f}:H_1(N-L)\to\Z_n$, where $\Z_n$ is the cyclic group of
order $n$. Moreover, two $n$-fold cyclic coverings $f':M'\to N$
and $f'':M''\to N$, branched over $L$, with associated monodromies
$\o_{f'}: H_1(N-L)\to\Z_n$ and $\o_{f''}: H_1(N-L)\to\Z_n$
respectively, are equivalent if and only if there exists
$u\in\Z_n$, with $\gcd(u,n)=1$, such that $\o_{f''}=u\o_{f'}$,
where $u\o_{f'}$ is the multiplication of $\o_{f'}$ by $u$.

Our goal is to obtain cyclic branched coverings with cyclically
presented fundamental groups. In order to achieve this, we will
select branched cyclic coverings of $(1,1)$-knots of ``special
type'', naturally generalizing the case of knots and links in
$\S^3$.

An $n$-fold cyclic covering $f$ of $L(p,q)$ branched over a
$(1,1)$-knot $K\subset L(p,q)$ will be called {\it
strongly-cyclic\/} (and denoted by $C_n(K)$) if the branching
index of $K$ is $n$. This means that the fiber $f^{-1}(x)$ of each
point $x\in K$ contains a single point. In this case the homology
class of a meridian loop $m$ around $K$ is mapped by $\o_{f}$ in a
generator of $\Z_n$ (up to equivalence we can always suppose
$\o_{f}(m)=1$). In the case of $(1,1)$-knots in $\S^3$,
strongly-cyclic branched coverings and cyclic branched coverings
are equivalent notions.

This type of covering appears to be a suitable tool for producing
3-manifolds with fundamental group admitting cyclic presentation,
as shown in the next section. For example, it is easy to see that
all Dunwoody manifolds are coverings of this type.

\medskip

As is well known, an $n$-fold cyclic branched covering of a knot
$K$ in $\S^3$ always exists and is unique (up to equivalence) for
all $n>1$, since $H_1(\S^3-K)\cong\Z$ and the homology class $m$
of a meridian loop around the knot is mapped by $\o_{f}$ in a
generator of $\Z_n$, which can be chosen equal to one, up to
equivalence. Obviously, this property is no longer true when the
ambient space is not a 3-sphere. For a $(1,1)$-knot $K\subset
L(p,q)$ the homology of $L(p,q)-K$ has the structure described in
Corollary \ref{homology}. As a consequence, the $n$-fold
strongly-cyclic branched covering of a $(1,1)$-knot not belonging
in $\S^3$ may not exist, and when it exists it may not be unique.

The next theorem gives necessary and sufficient conditions for the
existence of these coverings, establishing how many of these
coverings exist, up to equivalence.

\begin{theorem} \label{existence}
Let $K_\f$ be a $(1,1)$-knot in $L(p,q)$. Then $K_\f$ admits an
$n$-fold strongly-cyclic branched covering if and only if $d$
divides $q''$, where $d=\gcd (p,n)$ and $q''$ is uniquely
determined by the homology relation $ \f(\b)=p\a+q'\b+q''\bb$. In
this case there exists exactly $d$  of such coverings, up to
equivalence.
\end{theorem}
\begin{proof} By definition, the existence of the $n$-fold strongly-cyclic branched covering $f$
means that $\o_f(\bb)=1$, up to equivalence.  From Corollary \ref{homology}
we have $H_1(L(p,q)-K_{\f})=<\a,\bb\,|\,p\a+q''\bb>$, where $q''$
is uniquely determined by the homology relation
$\f(\b)=p\a+q'\b+q''\bb$. So, the covering exists if and only if
$\o_f(p\a+q''\bb)\equiv 0$ mod $n$ or, in other words, if and only
if there exists an element $x\in\Z_n$ such that $px+q''\equiv 0$ mod $n$.
This equation is solvable if and only if $d$
divides $q''$, where $d=\gcd (p,n)$ and in this case it has
exactly $d$ solutions. Since two different solutions give
non-equivalent coverings, the statement is proved.
\end{proof}

\medskip

In particular, for $(1,1)$-knots in $\S^3$ we have $p=\pm 1$ and
the existence and uniqueness of the $n$-fold cyclic
branched covering immediately follows.

\bigskip

\noindent {\bf Examples:}
\begin{itemize}
\item[(i)] let $\f=\d_{\a}^2\d_{\bb}\d_{\a}^{-4}$, then $K_{\f}$
is a $(1,1)$-knot in $L(6,5)$. Applying (\ref{formuleMCG}), we
obtain $\f (\lb)=(\la^2\lbb\la^{-1})^4\la^2\lb$. So we have
$$H_1(L(6,5)-K_{\f})=<\a ,\bb\,|\,6\a+4\bb >\cong \Z\oplus\Z_2.$$
In this case no $6$-fold strongly-cyclic branched coverings of
$K_{\f}$  exist, because $d=6$ does not divide $q''=4$. Observe
that there exists a 6-fold cyclic branched covering of $K_{\f}$:
take for example $\o_f(\a)=1$ and $\o_f(\bb)=3$; but obviously it
is not strongly-cyclic because the index around the knot is two;
\item[(ii)] let $\f=\d_{\a}^{-2}\d_{\bb}^{-2}\d_{\a}^{-2}$, then
$K_{\f}$ is a $(1,1)$-knot in $L(4,1)$. We have
$\f (\lb)=((\lbb^{-1}\la^2)^2\la^{-1})^2\la^{-2}\lb$ and therefore
$$H_1(L(4,1)-K_{\f})=<\a ,\bb\,|\,4\a-4\bb >\cong
\Z\oplus\Z_4.$$ In this case there are exactly four non-equivalent
$4$-fold strongly-cyclic branched coverings of $K_{\f}$, depending
on the choice of $\o_{f}(\a)\in\Z_4$.
\end{itemize}

Non-equivalent $n$-fold strongly-cyclic branched coverings of the
same $(1,1)$-knot may be effectively non-homeomorphic, as shown in
the Remark in the next section.

\section{Connections with cyclic presentations of  groups}

In this section we study the connections between strongly-cyclic
branched coverings of $(1,1)$-knots and cyclic presentations of
their fundamental groups. The following result has been announced
in \cite{Mu}.

\begin{theorem} \label{Main Theorem}
Every $n$-fold strongly-cyclic branched covering of a $(1,1)$-knot
admits a cyclic presentation for the fundamental group induced by
a Heegaard splitting of genus $n$.
\end{theorem}
\begin{proof}
Let $f:(M,f^{-1}(K))\to (L(p,q),K)=(T,A)\cup_{\ps}(T',A')$ be an
$n$-fold strongly-cyclic branched covering of the $(1,1)$-knot
$K$. Then $Y_n=f^{-1}(T)$ and $Y_n'=f^{-1}(T')$ are both
handlebodies of genus $n$. Moreover, $f^{-1}(A)$ and $f^{-1}(A')$
are both properly embedded trivial arcs in $Y_n$ and $Y_n'$
respectively. We get a genus $n$ Heegaard splitting
$(M,f^{-1}(K))=(Y_n,f^{-1}(A))\cup_{\F}(Y'_n,f^{-1}(A'))$, where
$\F:\partial Y_n'\to\partial Y_n$ is the lifting of $\ps$ with
respect to $f$. Let $m$ be a meridian loop around $A$ and
$\wa\subset T-A$ be a generator of $\pi_1(T,*)$, such that
$\o_{f}(\wa)=0$. It exists: take a generator $\eta\subset T-A$ of
$\pi_1(T,*)$; if $\o_{f}(\eta)=k$ then choose a loop $\wa$
homotopic to $\eta m^{-k}$. The set $f^{-1}(\wa)$ has exactly $n$
components $\ba_1,\ldots\ba_n$ and it is a set of generators for
$\pi_1(Y_n)$. A generator $\Psi$ of the group of covering
transformations cyclically permutes these components. Let $E'$ be
a meridian disk for the torus $T'$ such that $E'\cap
A'=\emptyset$, then $f^{-1}(E')$ is a system of meridian disks
$\{\tilde E_1',\ldots,\tilde E_n'\}$ for the handlebody $Y_n'$,
and they are cyclically permuted by $\Psi$. The curves
$\F(\partial\tilde E_1'),\ldots,\F(\partial\tilde E_n')$ give the
relators for the presentation of $\pi_1(M)$ induced by the
Heegaard splitting. Since both generator and relator curves are
cyclically permuted by $\Psi$, we get the statement.
\end{proof}

Since 2-bridge knots and torus knots are $(1,1)$-knots in $\S^3$
we have the following consequence:

\begin{corollary} \label{twobridge-torus}
Every $n$-fold cyclic branched covering of a 2-bridge knot and of
a torus knot admits a geometric cyclic presentation for the
fundamental group with $n$ generators.
\end{corollary}

We remark that the previous result for 2-bridge knots has also
been obtained (by another technique) in \cite{GM} and  the result
for torus knots largely generalizes a result obtained in
\cite{CHK1} and \cite{CHK2}.

\medskip

Let $K$ be a $(1,1)$-knot in $L(p,q)$ and let $\f$ be an element
of $PMCG(T_2)$ such that $K=K_{\f}$. Then
$\p_1(L(p,q)-K,*)=<\la,\lbb\,|\,r(\la,\lbb)>$, as shown in
Proposition \ref{fundamental}. Now, let $\o_{f}$ be the monodromy
of an n-fold strongly-cyclic branched covering of $K$. Following
the proof of Theorem \ref{Main Theorem}, we choose a new generator
$\wa =\la\lbb^{-\o_{f}(\bar\a)}$ and we get
$\p_1(L(p,q)-K,*)=<\wa,\lbb\,|\,\bar r(\wa,\lbb)>$, with $\bar
r(\wa,\lbb)= r(\wa\lbb^{\o_{f}(\bar\a)},\lbb)$. We have $\bar
r(\wa,\lbb)=\wa^{\ee_1}\lbb^{\dd_1}\cdots\wa^{\ee_s}\lbb^{\dd_s}$
for some $\ee_1,\ldots,\ee_s,\dd_1,\ldots,\dd_s\in\Z$.

\begin{theorem} \label{word} With
the assumptions listed in this section, let $\bar
r(\wa,\lbb)=\wa^{\ee_1}\lbb^{\dd_1}\cdots\wa^{\ee_s}\lbb^{\dd_s}$.
Then the fundamental group of the $n$-fold strongly-cyclic
branched covering of $K$, with monodromy $\o_f$, admits the cyclic
presentation $G_n(w)$, with:
$$w=x_{i_1}^{\ee_1}\cdots x_{i_s}^{\ee_s}$$
(subscripts mod $n$), where $i_k\equiv 1+\sum_{j=1}^{k-1}\dd_j
\mod n$, for $k=1,\ldots,s$.
\end{theorem}
\begin{proof}
From the proof of Theorem \ref{Main Theorem},  the fundamental
group of $C_n(K)$ has generators corresponding to the components
$\ba_1,\ldots,\ba_n$ of the lifting of $\wa$. The relators are
(homotopic to) $\F(\partial \tilde E'_1),\ldots, \F(\partial
\tilde E'_n)$, where $\F$ is the lifting of $\ps$ with respect to
$f$ and each $\partial \tilde E'_i$ is a component of the lifting
of a meridian disc $E'$ of $T'$ such that $T'\cap A'=\emptyset$.
So we can choose $E'$ such that $\partial E'=\tau(\b)$ and
therefore the relators are (homotopic to) the components of the
lifting of $r(\la,\lbb)$, or equivalently $\bar r(\wa,\lbb)$,
which arise from the relation $\lb'=\ps(\lb')$. Now, since
$\o_{f}(\bb)=1$, a factor $\lbb^k$ lifts to a path connecting the
point of $f^{-1}(*)$ in the sheet $i$ with the corresponding point
in the sheet $i+k$ (mod $n$). The result immediately follows.
\end{proof}

\medskip
We resume the algorithm for finding the word $w$ of the cyclic
presentation of $\pi_1(C_n(K))$, defined by the monodromy
$\o_{f}$, starting from an element $\f\in PMCG(T_2)$ such that
$K=K_{\f}$:
\begin{itemize}
\item[(i)] use (\ref{formuleMCG}) to calculate the homotopy class
$\f(\lb)$;
\item[(ii)] obtain $r(\la,\lbb)$ by erasing all $\lb$'s from
$\f(\lb)$;
\item[(iii)] compute $\bar r(\wa,\lbb)$ by replacing $\la$ with $\wa\lbb^{\o_{f}(\la)}$ in
$r(\la,\lbb)$;
\item[(iiii)] get $w$ by applying Theorem \ref{word}.
\end{itemize}

As an application we give an explicit cyclic presentation for
the fundamental group of the strongly-cyclic branched
coverings of the torus knots $\mathbf{t}(k,hk\pm 1)$, with $h,k>0$.

\begin{proposition} \label{proposition 1}
 The fundamental group of the $n$-fold cyclic branched covering
of the torus knot $\mathbf{t}(k,hk+1)$, with $h,k>0$, admits the
cyclic presentation $G_n(w)$, where
\begin{eqnarray*}
w&=&(\prod_{j=0}^{h(k-1)}x_{1-jk})(\prod_{i=0}^{k-2}(\prod_{l=1}^{h}x_{2+i-(h(k-1-i)+1-l)k}^{-1})),
\end{eqnarray*}
(subscripts mod $n$).\\
The fundamental group of the $n$-fold cyclic branched covering of
the torus knot $\mathbf{t}(k,hk-1)$, with $h,k>0$, admits the
cyclic presentation $G_n(w)$, where
\begin{eqnarray*}
w&=&(\prod_{j=1}^{h(k-1)-1}x_{1+jk}^{-1})(\prod_{i=0}^{k-2}(\prod_{l=0}^{h-1}x_{2+i+(h(k-1-i)-1-l)k})),
\end{eqnarray*}
(subscripts mod $n$).
\end{proposition}
\begin{proof} First we consider the case $\mathbf{t}(k,hk+1)$. From Proposition \ref{proposition 0} we have $\mathbf{t}(k,hk+1)=K_{\f}$ with
$\f=\d_{\a}^{-h}\d_{\bb}^{-k}\d_{\b}^{k+1}\d_{\a}$ and
$\p_1(\S^3-\mathbf{t}(k,hk+1))=<\la,\lbb\,|\,r(\la,\lbb)>,$ with
$r(\la,\lbb)=\la^{hk+1}(\lbb\la^{-h})^k$ (see the proof of the
proposition). We obtain
$H_1(\S^3-\mathbf{t}(k,hk+1))=<\a,\bb\,|\,\a+k\bb>$. Since
$\o_f(\lbb)=1$,  then $\o_f(\la)=-k$ and $\wa=\la\lbb^k$. We get
$\p_1(\S^3-\mathbf{t}(k,hk+1))=<\wa,\lbb\,|\,\bar r(\wa,\lbb)>,$
with $\bar
r(\wa,\lbb)=(\wa\lbb^{-k})^{1+hk}(\lbb(\lbb^k\wa^{-1})^h)^k$. So
the fundamental group of the $n$-fold strongly-cyclic branched
covering $\p_1(C_n(\mathbf{t}(k,h)))$ admits a cyclic presentation
$G_n(w)$, where
$w=(\prod_{j=0}^{hk}x_{1-jk})(\prod_{i=0}^{k-1}(\prod_{l=0}^{h-1}x_{2+i-(h(k-i)-l)k}^{-1}))$
(subscripts mod $n$). Since the last $h$ letters are the inverse
of the first $h$ (in the opposite order), we can remove them from
the word and, shifting all the indexes of the letters by $hk$, the
statement holds. The case $\mathbf{t}(k,hk-1)$ can be obtained in
a perfectly analogous way.
\end{proof}

\medskip

In particular, for $h=1$ and $k=2$  we have that
$\mathbf{t}(k,hk+1)=\mathbf{t}(2,3)$ is the trefoil knot and the
fundamental group of its $n$-fold cyclic branched covering is
$G_n(x_{3}x_{1}x_{2}^{-1})$, which is clearly isomorphic to the
Sieradski group $S(n)$ described in Section 1.

The above proposition clarifies the group presentations discussed
in the main theorem of \cite{CHK2}.

\medskip

 \textit{Remark.} Using the above algorithm it is easy to
show that non-equivalent $n$-fold strongly-cyclic branched
coverings of the same $(1,1)$-knot may be effectively
non-homeomorphic. For example, the fundamental groups of the four
non-equivalent 4-fold strongly cyclic coverings of  Example (ii)
in Section 4 admit cyclic presentations $G_4(w_i)$, $i=0,1,2,3$,
defined by the words:
$$\begin{array}{lcl}
w_0=x_4^2x_3x_2^2x_1^{-1},&\  &w_1=x_4x_1^3x_2x_1^{-1},\\
w_2=x_4x_2x_3x_4x_2x_1^{-1},&\ &w_3=x_4x_3x_1x_3x_2x_1^{-1}\end{array}
$$
for $\o_{f}(\a)=0,1,2,3$ respectively.  The first homology groups
of the strongly-cyclic branched coverings $C_4(K_{\f})$ are:
$$ H_1(C_4(K_{\f}))=\begin{cases}
     \Z\oplus\Z_4& \text{if }\ \ \o_{f}(\a)=1,3, \\
    \Z_8\oplus\Z_8 & \text{if }\ \ \o_{f}(\a)=0,2.
  \end{cases}$$
So at least two of these coverings are not homeomorphic.

\medskip

There is a strict relation between the Alexander polynomial
of a $(1,1)$-knot in $\S^3$ and the polynomial associated with the cyclic
presentation constructed according to Theorem \ref{word}.

\begin{proposition} \label{Alexander}
Let $K\subset\S^3$ be a $(1,1)$-knot. If $\Delta_K(t)$ is the
Alexander polynomial of $K$ and $f_w(t)$ is the polynomial
associated to the cyclic presentation of the $n$-fold cyclic
branched covering of $K$, obtained  by applying Theorem
\ref{word}, then $\Delta_K(t)=f_w(t)$, up to units of $\Z [t,
t^{-1}]$, when $n>\deg\Delta_K(t)$.
\end{proposition}
\begin{proof}
The statement follows from the arguments of Remarks 3 and 4 of
\cite{Mi}.
\end{proof}

\medskip

Since Dunwoody manifolds are strongly-cyclic branched coverings of
$(1,1)$-knots, the above result gives a positive answer to a
question posed in \cite{Du}. Moreover, it could be a possible
starting point to extend the notion of Alexander polynomial to all
$(1,1)$-knots.

\bigskip

\noindent \textit{Acknowledgements.} The authors  would like to
thank Andrei Vesnin for his helpful suggestions. Work performed
under the auspices of the G.N.S.A.G.A. of I.N.d.A.M. (Italy) and
the University of Bologna, funds for selected research topics.

\end{document}